\documentclass[12pt,a4paper,reqno]{article}
\usepackage{fleqn,amssymb,amscd,a4,enumerate,bm,amsthm,amsmath,float,multicol,}

\usepackage[a4paper,textwidth=5.25in,textheight=8.25in,centering]{geometry}
\theoremstyle{plain}
\newtheorem{theorem}{\bf Theorem}[section]

\theoremstyle{remark}

\newtheorem{proposition}{\bf Proposition}[section]

\usepackage{titlesec,soul}
\titleformat{\section}[block]{\large\bfseries}{\thesection.}{.5em}{}
\titlespacing*{\section}{0pt}{18pt}{15pt}
\date{}

\def\ge{\geqslant}
\begin{document}

\thispagestyle{empty} \noindent {\large\bf Cohomology algebra of the
orbit space of some free actions on spaces of cohomology type
$\bm{(a,b)}$}

\vskip1.5em
\noindent{\bf Hemant Kumar Singh}\\[.5em]
Department of Mathematics, University of Delhi\\
Delhi 110 007, India\\
Email: hksinghdu@gmail.com

\vskip1.5em
\noindent{\bf Tej Bahadur Singh}\\[.5em]
Department of Mathematics, University of Delhi\\
Delhi 110 007, India\\
Email: tej\_b\_singh@yahoo.co.in

\baselineskip23pt
\begin{description}
\item[\textbf{Abstract.}]

{\bf Let $X$ be a finitistic space with non-trivial cohomology
groups $H^{in} (X; \mathbb{Z}) \cong \mathbb{Z}$ with generators
$v_i$, where $i =0,1,2,3$. We say that $X$ has cohomology type
$(a,b)$ if $v_1^2 =av_2$ and $v_1 v_2 = bv_3$. In this note, we
determine the mod 2 cohomology ring of the orbit space $X/G$ of a
free action of $G=\mathbb{Z}_2$ on $X$, where both $a$ and $b$ are
even. In this case, we observed that there is no equivariant map
$\mathbb{S}^m\rightarrow X$ for $m>3n$, where $\mathbb{S}^m$ has the
antipodal action. Moreover, it is shown that $G$ can not act freely
on space $X$ which is of cohomology type $(a,b)$ where $a$ is odd
and $b$ is even. We also obtain the mod 2 cohomology ring of the
orbit space $X/G$ of free action of $G=\mathbb{S}^1$ on the space
$X$ of type $(0,b)$.}

\end{description}

\begin{description}
\item[\textbf{Key Words:}]

Space of type $(a,b)$, Free action,
Cohomology algebra, Spectral sequence
\end{description}

\begin{description}
\item[\textbf{2000 Mathematics Subject classification.}]
Primary 57S17; Secondary 57S25.
\end{description}

\thispagestyle{empty}

\vfill
\eject
\baselineskip23pt
\section{Introduction}

Let $G$ be a topological group acting continuously on a topological
space $X$. The set $\hat{x} =\{ gx|g \in G\}$ is called the orbit of
$x$. The set of all orbits $\hat{x}, x \in X$, is denoted by $X/G$
and assigned the quotient topology induced by the natural projection
$\pi : X \to X/G, x \to \hat{x}$. An action of $G$ on $X$ is said to
be free if for all $x \in X$ the isotopy subgroup $G_x =\{g \in
G|g(x) =x\}$ is the identity element $e$ of the group $G$. Let $X$
be a finitistic space with non-trivial cohomology groups $H^{in} (X;
\mathbb{Z}) \cong \mathbb{Z}$ with generators $v_i$, where $i
=0,1,2,3$. We say that $X$ has cohomology type $(a,b)$ if $v_1^2
=av_2$ and $v_1 v_2 = bv_3$. The spaces of type $(a,b)$ were
introduced by Toda [8]. Note that if $b \neq 0 (\bmod \,p)$, then
either $X$ has mod $p$ cohomology isomorphic to
$\mathbb{S}^n\times\mathbb{S}^{2n}$ product of spheres or a
projective space of dimension $3n$.  Throughout this paper, $H^*(X)$
will denote the \v{C}ech cohomology of space $X$.

Let $G =\mathbb{Z}_p$, $p$ an odd prime, act freely on a space $X$
of cohomology type $(a,b)$. If $a=0\,(\bmod\, p)$ then the possible
ring structure of the orbit space $X/G$ has been investigated in
[3],[4], and recently, we have determined the mod 2 cohomology ring
of $X/G$ when $G=\mathbb{Z}_2$ act freely on a space $X$ which is of
cohomology type $(a,b)$, where both $a$ and $b$ are odd [6].

It is easy to see that the space
$X=\mathbb{S}^4\cup_{\mathbb{S}^1}\mathbb{S}^6$ is of type $(0,0)$
and admits a free involution. In this note, we determine the mod 2
cohomology algebra of orbit space of free $G =\mathbb{Z}_2$ action
on spaces of cohomology type $(a,b)$, where both $a$ and $b$ are
even.  In this case, we observed that there is no equivariant map
$\mathbb{S}^m\rightarrow X$ for $m>3n$, where $\mathbb{S}^m$ has the
antipodal action. Moreover, we show that $G$ can not act freely on
spaces of cohomology type $(a,b)$, where $a$ is odd and $b$ is even.
We have also obtain the $\bmod\, 2$ cohomology ring of the orbit
space of a free $\mathbb{S}^1$-action on the spaces of type $(0,b)$.
We recall that a paracompact Hausdorff space is finitistic if every
open covering has a finite-dimensional refinement.
\section{Preliminaries}

In this section, we recall some known facts which will be used in
the proof of our theorems. Given a $G$-space $X$, there is an
associated fibration $X \overset{i}{\rightarrow} X_G
\overset{\pi}{\rightarrow} B_G$, and a map $\eta : X_G \to X/G$,
where $X_G =(E_G \times X)/G$ and $E_G \to B_G$ is the universal
$G$-bundle. When $G$ acts freely on $X, \eta: X_G \to X/G$ is
homotopy equivalence so the cohomology rings $H^*(X_G)$ and
$H^*(X/G)$ are isomorphic. To compute $H^*(X_G)$, we exploit the
Leray-Serre spectral sequence of the map $\pi : X_G \to B_G$. The
$E_2$-term of this spectral sequence is given by
\begin{eqnarray*}
E_2^{k,l} \cong H^k (B_G; \mathcal{H}^l (X))
\end{eqnarray*}
(Where $\mathcal{H}^{l}$ is a locally constant sheaf with stalk
$H^l(X)$ and group $G$) and it converges to $H^*(X_G)$, as
an algebra. The cup product in $E_{r+1}$ is induced
from that in $E_r$ via the isomorphism $E_{r+1} \cong H^*(E_r)$.
When $\pi_1(B_G)$ operates trivially on $H^*(X)$, the system
of local coefficients is simple (constant) so that
\begin{eqnarray*}
E_2^{k, l} \cong H^k (B_G) \otimes H^{l} (X).
\end{eqnarray*}

In this case, the restriction of the product structure in the
spectral sequence of the subalgebras $E_2^{*,0}$ and $E_2^{0,*}$
gives the cup products on $H^*(B_G)$ and $H^*(X)$ respectively. The
edge homomorphisms
\begin{eqnarray*}
H^p (B_G ) =E_2^{p,0} \to E_3^{p,0}\to \ldots \to E_{p+1}^{p,0}
=E_{\infty}^{p,0} \subseteq H^p (X_G),\mbox{and}
\end{eqnarray*}
\begin{eqnarray*}
H^q (X_G)\to E_{\infty}^{0,q} \subset \ldots \subset E_2^{0,q} =H^q (X)
\end{eqnarray*}
are the homomorphisms
\begin{eqnarray*}
\pi^{*} : H^p (B_G) \to H^q (X_G), \mbox{and}
\end{eqnarray*}
\begin{eqnarray*}
i^*: H^q (X_G)\to H^q(X)
\end{eqnarray*}
respectively. These results about spectral sequences can be found
[5]. The following facts are well known.

Recall that if $G =\mathbb{Z}_2$ or $\mathbb{S}^1$ then
\begin{eqnarray*}
H^*(B_G: \mathbb{Z}_2) =\begin{cases} \mathbb{Z}_2[t]& \det t =1, G =\mathbb{Z}_2\\[4pt]
\mathbb{Z}_2 [t]  &\deg t =2, G = \mathbb{S}^1\end{cases}
\end{eqnarray*}
\begin{proposition} Suppose $G =\mathbb{Z}_2$ or $\mathbb{S}^1$ acts freely on a
finitistic space $X$ with $H^j (X; \mathbb{Z}_2) =0$ for all $j >n$.
Then $H^j (X_G ;\mathbb{Z}_2 )=0$ for all $j>n$.
\end{proposition}

For $G=\mathbb{Z}_2$ see Bredon [1], and for $G=\mathbb{S}^1$ see
[7].

\section{$\bm{\mathbb{Z}_2}$ actions on spaces of type $\bm{(a,b)}$}\medskip

For the action of $G =\mathbb{Z}_2$ on $X$ of cohomology type
$(a,b)$, we obtain

\begin{theorem} Let $X$ be a finitistic space of cohomology type $(a,b)$ where $a$ is odd and $b$ is even. Then $G
=\mathbb{Z}_2$ can not act freely on $X$.
\end{theorem}

\begin{proof} By the Universal coefficient theorem, $H^{in} (X; \mathbb{Z}_2) =\mathbb{Z}_2$,
for $i =0,1,2,3$. Let $v_i\epsilon H^{in} (X ;\mathbb{Z}_2)$, $i
=1,2,3$ be generators. Then we have the relation $v_1^2 =v_2$ and
$v_1 v_2 =0$. As $\pi_1 (B_G) =\mathbb{Z}_2$ acts trivally on
$H^*(X)$, the fibration $X \overset{i}{\rightarrow} X_G
\overset{\pi}{\rightarrow} B_G$ has a simple system of local
coefficients on $B_G$. Therefore, the spectral sequence has
\begin{eqnarray*}
E_2^{k,l } \cong H^k (B_G) \otimes H^{l } (X)
\end{eqnarray*}
If $G =\mathbb{Z}_2$ acts freely on $X$ then the Lerray-Serre
spectral sequence of the map $\pi : X_G \to B_G$ must not collapse
at $E_2$-term. So the differentials
\begin{eqnarray*}
d_r :E_r^{k,l} \to E_r^{k+r, l -r -1}
\end{eqnarray*}
can not be trivial for all $r$. The non-trivial differentials are
possible only for $r =n+1, 2n+1$ and $3n +1$.  It is obvious that
$E_2^{k,l} =\mathbb{Z}_2$ for all $k$ when $l = 0, n, 2n, 3n$;  and
$E_2^{k,l} =0$, otherwise. Consequently, $E_2^{k,l } =E_n^{k,l }$
for $k$ and $l $. If $d_{n+1} (1\otimes  v_1) =t^{n+1} \otimes 1$.
Then $d_{n+1}(1\otimes v_2)=0$ so that $0 = d_{n+1} ((1 \otimes v_1)
(1 \otimes v_2))= t^{n+1} \otimes v_2$. This contradiction forces
that $d_{n+1}(1\otimes v_1)$ must be trivial. If
\begin{eqnarray*}
d_{n+1} :E_{n+1}^{0,2n} \to E_{n+1}^{n+1,n}
\end{eqnarray*}
is non trivial, then $d_{n+1} (1\otimes v_2) =t^{n+1} \otimes v_1$.
So $0 =d_{n+1}((1\otimes v_1)(1\otimes v_2)) =t^{n+1} \otimes v_1^2
=t^{n+1} \otimes v_2$, a contradiction. Therefore, $d_{n+1}(1\otimes
v_2)$ is also trivial. Now, if $d_{n+1}(1\otimes v_3)\neq 0$, then
two lines in the spectral sequence survive to infinity, and this
contradicts Proposition 2.1. And, if $d_{n+1}(1\otimes v_3)= 0$,
then $d_{2n+1}(1\otimes v_2)= 0$ and it follows that at least two
lines in the spectral sequence survive to infinity, again
contradicting Proposition 2.1. This completes the proof.
\end{proof}

For the spaces of type $(a,b)$ where both $a$ and $b$ are even, we
obtain
\begin{theorem} Let $G =\mathbb{Z}_2$ act freely on a finitistic space of cohomology type
$(a,b)$ where both $a$ and $b$ are even. Then as a graded
Commutative algebra
\begin{eqnarray*}
H^*(X/G ; \mathbb{Z}_2) =\mathbb{Z}_2 [x,z] / \langle x^{3n+1}, z^2, zx^{n+1} \rangle
\end{eqnarray*}
where $\deg x =1$ and $\deg z =n$.
\end{theorem}

\begin{proof} As in Theorem 3.1
\begin{eqnarray*}
E_2^{k,l}\cong H^{k}(B_G) \otimes H^{l} (X).
\end{eqnarray*}
Since $G=\mathbb{Z}_2$ acts freely on $X$ so the spectral sequence
of map $\pi : B_G\rightarrow X_G$ can neither degenerate nor any
line can survive to infinity. It is easy to see that
\begin{eqnarray*}
d_{n+1}(1\otimes v_1)=d_{n+1}(1\otimes v_3)=0 \, \text{and}\,
d_{n+1}(1\otimes v_2)\neq 0
\end{eqnarray*}

Therefore, we get $E_{n+2}^{k,l } =\mathbb{Z}_2$ for all $k$ if $l
=0, 3n$; for $k =0,1,2,\ldots, n$ if $l =n$ and $E_n^{k,l } =0$,
otherwise. Now, it is obvious that $d_{2n+1}=0$ and
$d_{3n+1}(1\otimes v_3)\neq 0$. Thus $E_{3n+2}^{k,l } =\mathbb{Z}_2$
for $k\leq 3n$ if $l =0$; $k \leq n$ if $l =n$. It follows that
\begin{eqnarray*}
H^j(X_G) =\begin{cases} \mathbb{Z}_2 \oplus \mathbb{Z}_2 & \text{for} \ n \leq j \leq 2n\\[4pt]
\mathbb{Z}_2 & \text{for}  \ \ 0\leq j \leq n-1 \ \text{and} \ \ 2n+1\leq j \leq 3n\\[4pt]
0 & \text{otherwise}\end{cases}
\end{eqnarray*}
Now, we compute the multiplication in $H^*(X_G)$. The element $1
\otimes v_1 \in E_2^{0,n}$ is permanent cocycle so this determine $z
\in E_{\infty}^{0,n}$ such that $i^*(z) =v_1$ and $z^2 =0$. Let $x
=t\otimes 1 \in E_{\infty}^{1,0}$. Then $x^{3n+1} =0$ and $zx^{n+1}
=0$. If follows that
\begin{eqnarray*}
\text{Tot} \ \  E_{\infty}^{*,*} =\mathbb{Z}_2[x,z] / \langle x^{3n+1}, z^2, zx^{n+1}\rangle
\end{eqnarray*}
where $\deg x =1, \deg z =n$.

Thus $H^*(X_G) =\mathbb{Z}_2[x,z]/\langle x^{3n+1}, z^2,
zx^{n+1}\rangle$ where $\deg x =1$ and $\deg z=n$. Since the action
of $G$ on $X$ is free, $\pi : X_G \to X/G$ is a homotopy equivalence
and so induces a cohomology isomophism. This completes the proof.
\end{proof}

\section{$\bm{\mathbb{S}^1}$ actions on spaces of type $\bm{(a,b)}$}\medskip

For actions of $\mathbb{S}^1$ on spaces of type $(a,b)$, we have
following result
\begin{theorem} Let $G =\mathbb{S}^1$ act freely on a finitistic space $X$ of cohomology type
$(a,b)$. Then $a$ must be zero, and $H^*(X/G;\mathbb{Z}_2)$ is one
of the following graded commutative algebras:
\begin{eqnarray*}
(i)\mathbb{Z}_2[x,z]/\langle x^{\frac{3n+1}{2}}, z^2,
zx^{\frac{n+1}{2}}\rangle
\end{eqnarray*}
where $\deg x=2$, $\deg z =n$.
\begin{eqnarray*}
(ii) \mathbb{Z}_2[x,z]/\langle x^{\frac{n+1}{2}}, z^2\rangle
\end{eqnarray*}
where $\deg x=2$, $\deg z =2n$ and $b$ is odd.
\end{theorem}

\baselineskip24pt
\begin{proof} Note that if $n$ is even,
then $\mathbb{S}^1$ can not act freely on $X$. Let $v_i\,\epsilon
H^{in}(X;\mathbb{Z}_2)$, $i =1,2,3$ be genarators, where $n$ is odd.
So $v_1^2=0$. Since $\pi_1(B_G)$ is trivial, so the $E_2$-term of
the Leray-Serre spectral sequence is
\begin{eqnarray*}
E_2^{k,l } =H^k (B_G)\otimes H^{l } (X)
\end{eqnarray*}
Clearly $E_2^{k,l } =\mathbb{Z}_2$ for even $k$ when $l  =0, n, 2n,
3n$; and $E_2^{k,l } =0$, otherwise. Since the cohomology of $X_G$
vanishes in high degrees the spectral sequence of the fibration $X
\overset{i}{\rightarrow} X_G \overset{\pi}{\rightarrow} B_G$ does
not collapse at $E_2$-term.

First, assume that $d_{n+1}(1\otimes v_1) =0$. Then, by Proposition
2.1, $d_{n+1}(1\otimes v_2)=t^\frac{n+1}{2}\otimes v_1$ and
$d_{n+1}(1\otimes v_3)= 0$. Obviously, $d_{2n+1}=0$. Again, by
Proposition 2.1, $d_{3n+1}(1\otimes v_3)\neq 0$. Now, we have
$E_{3n+2}^{k,l } =\mathbb{Z}_2$ for even $k\leq 3n-1$ if $l =0$; for
even $k\leq n-1$ if $l  =n$ and $E_{3n+2}^{k, l } =0$, otherwise. It
is obvious that $E_{3n+2} =E_{\infty}$. Thus, when $n=1$,
$H^j(X_G)=\mathbb{Z}_2$ for $j=0,1,2$ and $H^j(X_G)=0$, otherwise.
And, when $n>2$\mathindent0em
\begin{eqnarray*}
H^j(X_G) =\begin{cases} 0 &\!\!\! \text{for} \ j =2i +1 (0\leq i \leq \frac{n-3}{2} \ \text{or} \ n \leq i \leq \frac{3n-3}{2}) \text{ or } j\ge 3n\\
\mathbb{Z}_2 &\!\!\! \text{otherwise} \end{cases}
\end{eqnarray*}\mathindent3em

 Choose $x \in H^2(X_G)$ such that $\pi^*(t) =x$. Then $x^{\frac{3n+1}{2}} =0$ and the
multiplication by $x$
\begin{eqnarray*}
x\cup (\cdot): E_{\infty}^{k,l }\to E_{\infty}^{k +2,l }
\end{eqnarray*}
is an isomorphism for $k < 3n-1$ if $l  =0$ and for $k < n-1$ if $l
=n$. Therefore multiplication by $x \in H^2(X_G)$
\begin{eqnarray*}
x\cup (\cdot): H^k(X_G)\to H^{k+2} (X_G)
\end{eqnarray*}
is an isomorphism for even $k < n-1$. Then element $1 \otimes v_1
\in E_2^{0,n}$ is a permanent cocyle so determine an element $z \in
E_{\infty}^{0,n}$ such that $i^*(z) =v_1$ and $z^2 =0$. Clearly,
$zx^{\frac{n+1}{2}} =0$. If follows that the total complex Tot
$E_{\infty}^{*,*}$ is graded algebra given by
\begin{eqnarray*}
\text{Tot} E_{\infty}^{*,*} =\mathbb{Z}_2[x,z]/\langle
x^{\frac{3n+1}{2}}, z^2, zx^{\frac{n+1}{2}}\rangle
\end{eqnarray*}
where $\deg x=2$, $\deg z =n$.

Therefore, we have
\begin{eqnarray*}
H^*(X_G) =\mathbb{Z}_2[x,z]/\langle x^{\frac{3n+1}{2}},z^2,
zx^{\frac{n+1}{2}}\rangle
\end{eqnarray*}
where $\deg x =2$, $\deg z =n$.

Next, consider the case $d_{n+1}(1\otimes v_1)
=t^{\frac{n+1}{2}}\otimes 1$. Then $d_{n+1}(1\otimes v_2)=0$. Now,
if $b$ is even then $v_1v_2=0$ implies that $0=d_{n+1}((1\otimes
v_1)(1\otimes v_2))=t^{\frac{n+1}{2}}\otimes v_2$, a contradiction.
So , $b$ must be odd. Thus we have $d_{n+1}(1\otimes
v_3)=t^{\frac{n+1}{2}}\otimes v_2$. It follows that $E_{\infty}^{k,l
}=\mathbb{Z}_2$ for $k=0,2,4\ldots n-1$ if $l=0,2n$;
$E_{\infty}^{k,l }=0$, otherwise. Thus
  \begin{eqnarray*}
H^j(X_G) =\begin{cases} \mathbb{Z}_2 &\!\!\! \text{for}\,j=2i (0\leq
i \leq \frac{n-1}{2}\, \text{or}\, n\leq i \leq \frac {3n-1}{2})\\ 0
&\!\!\! \text{otherwise}.\end{cases}
\end{eqnarray*}

 The element $1\otimes v_2\,\epsilon E^{0,2n}_2$ is a permanent cocyle so
determine an element $z\,\epsilon\,E_\infty^{0,n}$ such that
$i^*(z)=v_2$ and $z^2=0$. If $n=1$, we have
\begin{eqnarray*}
H^*(X_G) =\mathbb{Z}_2[z]/\langle z^2\rangle, \text{where $\deg
z=2$.}
\end{eqnarray*}
Now, assume that $n>2$. Put $\pi^*(t)=x$. Then $x \epsilon
\,H^2(X_G)$ is determined by $t\otimes 1 \epsilon E^{2,0}_\infty$
and $x^{\frac{n+1}{2}}=0$. If follows that the total complex Tot
$E_{\infty}^{*,*}$ is graded algebra given by
\begin{eqnarray*}
\text{Tot} E_{\infty}^{*,*} =\mathbb{Z}_2[x,z]/\langle
x^{\frac{n+1}{2}}, z^2\rangle
\end{eqnarray*}
where $\deg x=2$,$\deg z =2n$ and $b$ is odd.
Therefore, we have
\begin{eqnarray*}
H^*(X_G) =\mathbb{Z}_2[x,z]/\langle x^{\frac{n+1}{2}},z^2\rangle.
\end{eqnarray*}
This completes the proof.
\end{proof}

\section{Index of spaces of cohomology type $(a,b)$}\medskip
Let $G=\mathbb{Z}_2$ acts freely on a finitistic space $X$. We
recall that an equivariant map from a $G$-space $X$ to a $G$-space
$Y$ is a continuous map $\phi :X\rightarrow Y$ such that
$g\phi(x)=\phi g(x)$ for all $g\epsilon G, x\epsilon X$. The index
of $X$ is defined to be the largest integer $n$ such that there
exists an equivariant map $\mathbb{S}^n\rightarrow X$ relative to
antipodal action on $n$-sphere $\mathbb{S}^n$. And, the mod 2
cohomology index of $X$ is the largest integer $n$ such that
$x^n\neq 0$, where $x\epsilon \, H^1(X/G;\mathbb{Z}_2)$ is the euler
class of the bundle $\mathbb{Z}_2\hookrightarrow X\rightarrow
X/\mathbb{Z}_2$, by Conner and Floyd in [2]. It has been shown that
index of $X$ can not exceeds the mod 2 cohomology index of $X$,
$(4.5, ibid)$.

If $a$ is odd and $b$ is even then observe that $G=\mathbb{Z}_2$ can
not act freely on $X$. By Theorem 3.2, the mod 2 cohomology index of
$X$ is $3n$ when both $a$ and $b$ are even. If both $a$ and $b$ are
odd, then the mod 2 cohomology index of $X$ is 2 (see Theorem 1 in
[6]), and when $a$ is even and $b$ is odd, the mod 2 cohomology
index of $X$ is $n$ or $3n$ (see Theorem 2 in [4]).

By the above remarks, it follows that with the antipodal action on
$\mathbb{S}^m$, there exists no equivariant map
$\mathbb{S}^m\rightarrow X$ for (i) $m>3n$ when both $a$ and $b$ are
even, and (ii) $m>2$ when both $a$ and $b$ are odd.

\end{document}